\documentclass[11pt, a4paper, twoside]{amsart}
\usepackage{amsmath}
\usepackage{amssymb}
\usepackage{amsfonts}
\usepackage{a4}
\usepackage[latin1]{inputenc}

\numberwithin{equation}{section}

\newtheorem{theorem}{Theorem}[section]
\newtheorem{lemma}[theorem]{Lemma}

\theoremstyle{definition} 

\theoremstyle{plain} 

\newtheorem{fact}[theorem]{Fact}
\newtheorem{conclusion}[theorem]{Conclusion}

\newtheorem{mainlemma}[theorem]{Main Lemma}

\newtheorem*{maintheorem*}{Main Theorem} 
\newtheorem*{conjecture*}{Conjecture} 
\newtheorem{definition}[theorem]{Definition}

\newtheorem{remark}[theorem]{Remark}

\theoremstyle{remark}  

\newcommand{\nc}{\newcommand}
\nc{\nothing}[1]{}

\nothing{   
\nc{\dom}{{\rm dom}}
\nc{\card}{{\rm card}}
\nc{\lh}{{\rm lh}}
\nc{\lgg}{{\rm lg}}
\nc{\rge}{\mbox{\rm range}}
\nc{\cf}{{\rm cf}}
\nc{\nex}{\mbox{\rm next}}
\nc{\uhr}{\restriction}
\nc{\supt}{{\rm supt}}
\nc{\supp}{{\rm supp}}
\nc{\Lim}{{\rm Lim}}
\nc{\Leb}{{\rm Leb}}
\nc{\modd}{{\rm mod}}
\nc{\RO}{{\rm RO}}
\nc{\prob}{{\rm Prob}}
}

\nc{\On}{{\rm On}}

\nc{\nco}{\DeclareMathOperator}

\nco{\order}{o} 
\nco{\ppower}{pp} 
\nco{\pcf}{pcf} 
\nco{\tcf}{tcf} 
\nco{\tlim}{tlim} 
\nco{\limtext}{lim} 
\nco{\prodt}{{\textstyle \prod}} 
\nco{\symdiff}{\triangle}
\nco{\dom}{dom}
\nco{\card}{card}
\nco{\lh}{lh}
\nco{\lgg}{lg}
\nco{\rge}{range}
\nco{\otp}{otp}
\nco{\trunk}{tr}
\nco{\cf}{cf}
\nco{\nex}{next}
\nc{\uhr}{\restriction}
\nco{\supt}{supt}
\nco{\supp}{supp}
\nco{\Lim}{Lim}
\nco{\Leb}{Leb}
\nco{\modd}{mod}
\nco{\invariant}{inv}
\nco{\id}{id}
\nco{\RO}{RO}
\nco{\Dp}{Dp} 
\nco{\pss}{ps}
\nco{\acc}{acc}

\nc{\potom}{\ensuremath{{\cal P}(\omega)}}
\nc{\potinf}{\ensuremath{[\omega]^\omega}}
\nc{\pfin}{\ensuremath{{\cal P}(\omega)/{\rm fin}}}

\nc{\potfin}{\ensuremath{[\omega]^{<\omega}}}
\nc{\inn}{\ensuremath{{\omega^{\uparrow \omega}}}}
\nc{\hoch}{^{<\omega}}
\nc{\hocho}{^{\omega}}
\nc{\tree}[1]{{[} #1 {]}_0}
\nc{\tre}[2]{ {#1}_{#2}}
   
\nc{\prooff}[1]{{\bf Proof} of #1:}
\nc{\proofend}{\makebox{} \hfill ${\bf \square}$ \\}
\nc{\proofendof}[1]{\makebox{} \hfill $\boldmath{\square}_{\rm #1}$ \\}
\nc{\beq}{\begin{eqnarray*}}
\nc{\eeq}{\end{eqnarray*}}
\nc{\bde}{\begin{list}}
\nc{\ede}{\end{list}}

\newenvironment{myrules}
{\begin{list}{}
{
 \setlength{\leftmargin}{0.5in}
 \setlength{\labelwidth}{1cm}
 \setlength{\labelsep}{0.2in}
 \setlength{\parsep}{0.5ex plus 0.2ex minus 0.1 ex}
 \setlength{\itemsep}{0.3ex plus 0.2 ex minus 0ex}
}}{\end{list}}

{\end{list}}

\newcounter{subalph}
{\end{list}}

\newcommand{\greek}[1]{\ifthenelse{\value{#1}=1}{\mbox{$\alpha$}}%
  {\ifthenelse{\value{#1}=2}{\mbox{$\beta$}}{%
   \ifthenelse{\value{#1}=3}{\mbox{$\gamma$}}{%
   \ifthenelse{\value{#1}=4}{\mbox{$\delta$}}{%
   \ifthenelse{\value{#1}=5}{\mbox{$\varepsilon$}}{%
   \ifthenelse{\value{#1}=6}{\mbox{$\zeta$}}{%
   \ifthenelse{\value{#1}=7}{\mbox{$\eta$}}{%
   \ifthenelse{\value{#1}=8}{\mbox{$\theta$}}{%
   \ifthenelse{\value{#1}=9}{\mbox{$\iota$}}{%
   \ifthenelse{\value{#1}=10}{\mbox{$\kappa$}}{%
   \ifthenelse{\value{#1}=11}{\mbox{$\lambda$}}{%
   \ifthenelse{\value{#1}=12}{\mbox{$\mu$}}{%
   \ifthenelse{\value{#1}=13}{\mbox{$\nu$}}{%
   \ifthenelse{\value{#1}=14}{\mbox{$\xi$}}{%
   \ifthenelse{\value{#1}=15}{\mbox{$\rm o$}}{%
   \ifthenelse{\value{#1}=16}{\mbox{$\pi$}}{%
   \ifthenelse{\value{#1}=17}{\mbox{$\varrho$}}{%
   \ifthenelse{\value{#1}=18}{\mbox{$\sigma$}}{%
   \ifthenelse{\value{#1}=19}{\mbox{$\tau$}}{%
   \ifthenelse{\value{#1}=20}{\mbox{$\upsilon$}}{%
   \ifthenelse{\value{#1}=21}{\mbox{$\varphi$}}{%
   \ifthenelse{\value{#1}=22}{\mbox{$\chi$}}{%
   \ifthenelse{\value{#1}=23}{\mbox{$\psi$}}{\mbox{$\omega$}%
  }}}}}}}}}}}}}}}}}}}}}}}}

\newcounter{subgreek}
{\end{list}}

\newcounter{subarabic}
{\end{list}}

\newcounter{subroman}
{\end{list}}

\newcount\skewfactor

\def\mathunderaccent#1#2 {\let\theaccent#1\skewfactor#2
\mathpalette\putaccentunder}
\def\putaccentunder#1#2{\oalign{$#1#2$\crcr\hidewidth
\vbox to.2ex{\hbox{$#1\skew\skewfactor\theaccent{}$}\vss}\hidewidth}}
\def\name{\mathunderaccent\tilde-3 }

\nc{\nname}{\name}

\nc{\even}{\ensuremath{\rm Even}}
\nc{\odd}{\ensuremath{\rm Odd}}

\nc{\al}{$\alpha$\  }
\nc{\om}{\omega}
\nc{\omm}{\ensuremath{\omega_1}}
\nc{\ep}{\varepsilon}
\nc{\tk}{\tilde{K}}
\nc{\concat}{\,\hat{} \,}   
\nc{\force}{\Vdash}
\nc{\fb}{f_{\bar{M}}}
\nc{\such}{\, : \,}   

\nc{\meager}{\ensuremath{{\cal M}}}
\nc{\lebesgue}{\ensuremath{{\cal N}}}
\nc{\nulll}{\ensuremath{{\cal N}}}
\nc{\ksigma}{\ensuremath{{\bf K}_\sigma}}
\nc{\ideal}{\ensuremath{{\cal I}}}
\nc{\ga}{\ensuremath{\frak a}}
\nc{\AAA}{{\cal A}}   
\nc{\gc}{\ensuremath{\frak c}}
\nc{\gs}{\ensuremath{\frak s}}
\nc{\gh}{\ensuremath{\frak h}}
\nc{\gd}{\ensuremath{\frak d}}
\nc{\gb}{\ensuremath{\frak b}}
\nc{\gro}{\ensuremath{\frak g}}
\nc{\gu}{\ensuremath{\frak u}} 
\nc{\gr}{\ensuremath{\frak r}} 
\nc{\gt}{\ensuremath{\frak t}}
\nc{\fff}{\ensuremath{\frak f}}
\nc{\gm}{\ensuremath{\mathfrak{mcf}}}
\nc{\gge}{\ensuremath{\mathfrak e}}
\nc{\cfupro}{\ensuremath{\cf(\upro)}}
\nc{\cfvpro}{\ensuremath{\cf(\vpro)}}
\nc{\gp}{\ensuremath{\frak p}}
\nc{\gk}{\ensuremath{\frak k}}

\nc{\add}[1]{\mbox{\ensuremath{{\rm add}(#1)}}}
\nc{\cov}[1]{\mbox{\ensuremath{{\rm cov}(#1)}}}
\nc{\unif}[1]{\mbox{\ensuremath{{\rm unif}(#1)}}}
\nc{\cof}[1]{{\mbox{\ensuremath{\rm cof}(#1)}}}

\nc{\addd}[2]{\mbox{\ensuremath{{\rm add}^{#1}(#2)}}}   
\nc{\covv}[2]{\mbox{\ensuremath{{\rm cov}^{#1}(#2)}}}   
\nc{\uniff}[2]{\mbox{\ensuremath{{\rm unif}^{#1}(#2)}}} 
\nc{\coff}[2]{{\mbox{\ensuremath{\rm cof}^{#1}(#2)}}}

\nc{\cd}{Cicho\'n's Diagram}
\nc{\COF}{\mbox{\bf Cof}}

\nc{\MA}{\mbox{\rm MA}}
\nc{\PFA}{\mbox{\rm PFA}}
\nc{\OCA}{\mbox{\rm OCA}}
\nc{\GCH}{\mbox{\rm GCH}}
\nc{\CH}{\mbox{\rm CH}}
\nc{\zfc}{\mbox{\rm ZFC}}
\nc{\sch}{\mbox{\rm SCH}} 
\nc{\ZF}{\mbox{\rm ZF}}
\nc{\NCF}{\mbox{\rm NCF}} 
\nc{\FD}{\mbox{\rm FD}}   

\nc{\fourG}{\mbox{\rm 4G}} 
\nc{\fourI}{\mbox{\rm 4I}}   

\nc{\Borelhood}{Borel measurability} 
\nc{\Pieinseins}{\mbox{${\bf \Pi}^1_1$}}
\nc{\seinseins}{\mbox{${\bf\Sigma}^1_1$}}
\nc{\seinszwei}{\mbox{${\bf\Sigma}^1_2$}}
\nc{\seinsdrei}{\mbox{${\bf\Sigma}^1_3$}}
\nc{\Deleinszwei}{\mbox{${\bf\Delta}^1_2$}}

\nc{\up}{\ensuremath{{\cal U}\mbox{\ensuremath{\rm -prod}}\,\omega}}
\nc{\upp}{\ensuremath{{\cal U}'\mbox{\ensuremath{\rm -prod}}\,\omega}}
\nc{\upro}{\ensuremath{{\cal U}\mbox{\ensuremath{\rm -prod}}\,\om}}
\nc{\fupro}{\ensuremath{f({\cal U})\mbox{\ensuremath{\rm -prod}}\,\om}}
\nc{\vpro}{\ensuremath{{\cal V}\mbox{\ensuremath{\rm -prod}}\,\om}}
\nc{\fpro}{\ensuremath{{\cal F}\mbox{\ensuremath{\rm -prod}}\,\om}}

\nc{\cff}[1]{{\text{cf}\,(#1)}}           
\nc{\cu}{\ensuremath{\cal U}}             
\nc{\ai}{\ensuremath{\forall^\infty}}     
\nc{\ei}{\ensuremath{\exists^\infty}}     
\nc{\ww}{\ensuremath{\omega^\omega}}      

\nc{\gw}{groupwise dense}

\nc{\kk}{car\-dinal cha\-rac\-teris\-tic}
\nc{\joker}{\ast}
\nc{\gtc}{Galois-Tukey connection} 

\nc{\av}[1]{{\rm Av}_{#1}}
\nc{\eps}{\varepsilon}
\nc{\n}{{\bf n}}                 
\nc{\m}{{\bf m}}

\nc{\marginparr}[1]{}
\nc{\footnoteee}{} 
\nc{\footnotee}{}  

\newcommand{\cal}{\mathcal}

\nc{\divs}{{c_0 \setminus \ell^1}}
\nc{\divser}{(\divs, \leq^*)/\thickapproy}
\nc{\bfin}{\RO(\pfin \setminus\{0\},\subseteq^*)}
\nc{\bdivser}{\RO(\divser)}
\nc{\inc}{{\rm INC}}
\nc{\com}{{\rm COM}}
\nc{\thickapproy}{\makebox{}\!\!\thickapprox}
\nc{\approy}{\makebox{}\!\!\approx}
\nc{\lessi}{\leqslant}
\nc{\gessi}{\geqslant}
\nc{\interior}[1]{{\rm int}(#1)}
\nc{\closure}[1]{{\rm cl}(#1)}
\nc{\Vo}{Vojt\'a\v{s}}

\nc{\precedeseq}{\leq^*} 
\nc{\precedes}{\prec}
\nc{\stronger}{\leqslant_{\bf P}}
\nc{\underlline}[1]{\hat{#1}}
\nc{\PO}{{\bf P}}
\nc{\charak}{\text{ch}}

\nc{\needed}{needed\ }
\nc{\neededc}{needed}
\nc{\Needed}{Needed\ }
\nc{\wneeded}{weakly needed\ }
\nc{\Wneeded}{Weakly needed\ }
\nc{\wneededc}{weakly needed}

\nc{\mup}{m_{\rm up}}
\nc{\mdn}{m_{\rm dn}}
\nco{\may}{may}
\nco{\aver}{av} 
\nco{\norm}{{\bf nor}} 
\nco{\val}{{\bf val}} 
\nco{\dis}{{\bf dis}} 
\nco{\basis}{basis}
\nco{\pos}{pos}

\nc{\err}{\mbox{err}}
\nc{\eee}{\mbox{e}}
\nco{\Expect}{Exp}

\begin{document}


\title{The splitting number can be smaller than the matrix chaos number}

\author{Heike Mildenberger and Saharon Shelah}
\thanks{The first author was supported by
a Minerva fellowship.}

\thanks{The second author's research 
was partially supported by the ``Israel Science
Foundation'', founded by the Israel Academy of Science and Humanities.
This is  the second author's work number 753}

\address{Heike Mildenberger,
Saharon Shelah,
Institute of Mathematics,
The Hebrew University of Jerusalem,
 Givat Ram,
91904 Jerusalem, Israel
}

\email{heike@math.huji.ac.il}
 
\email{shelah@math.huji.ac.il}

\begin{abstract}

Let $\chi$ be the minimum cardinal of
a subset of $2^\omega$ that cannot be made convergent
 by multiplication with  a single Toeplitz matrix.
By an application of creature forcing 
we show that $\gs < \chi$ is consistent.
We thus answer a question by \Vo.
We give two kinds of  models for the strict inequality.
The first is the combination of an
 $\aleph_2$-iteration of some proper
forcing with adding $\aleph_1$ random reals.
The second kind of models is got by adding $\delta$ 
random reals
to a model of $\MA_{<\kappa}$ for some $\delta \in [\aleph_1,\kappa)$.
It was a conjecture of Blass that
$\gs=\aleph_1 < \chi = \kappa$ holds in such a model.
For the analysis of the second model we again use the creature forcing
from the first model.
\end{abstract}



\maketitle


\nc{\LOC}{{\mathbb L}}

\setcounter{section}{-1}
\section{Introduction}
We consider products of $\omega\times \omega$ matrixes 
$A=(a_{i,j})_{i,j\in\omega}$ and  functions 
from $\omega$ to 2 or to some bounded interval of the reals.
The product $A \cdot f$ is defined
 as usual in linear algebra, i.e.,
$(A \cdot f)(i) = \sum_ {j \in \omega} a_{i,j} \cdot f(j)$.
We define
\[
A\lim f := \lim_{i\to\infty} \sum_{j=0}^{\infty} (a_{i,j}\cdot f(j)).
\]

\smallskip

Toeplitz (cf.\ \cite{Cooke}) showed:  $A\lim$ is an extension of the
ordinary limit iff $A$ is a regular  matrix\index{regular matrix}, i.e.\
iff $\exists m \; \forall i\; 
\sum_{j=0}^{\infty} |a_{i,j}| < m$ and
$\lim_{i \to \infty} \sum_{j=0}^{\infty} a_{i,j} = 1$ and
$\forall j \; \lim_{i \to \infty} a_{i,j} = 0$.
Regular matrices are also called Toeplitz\index{Toeplitz matrix} matrices.

\smallskip

We are interested whether for many $f$'s simultaneously there is one
$A$ such that all $A\lim f$ exist, and formulate our question in terms
of cardinal characteristics.

Let $\ell^\infty$ denote the set of bounded real
sequences, and let $\mathbb M$ denote the set of all Toeplitz
matrices.
\Vo\ \cite{Vojtas88}
defined for ${\mathbb A} \subseteq {\mathbb M}$
the chaos relations $\chi_{{\mathbb A}, \infty}$  and
their norms $\Vert \chi_{{\mathbb A}, \infty} \Vert$
\begin{eqnarray*}
\chi_{{\mathbb A},\infty} &=& \{(A,f) \such A \in {\mathbb A}
\;\wedge\; f \in  
\ell^\infty \;\wedge\; A \lim f \mbox{ does not exist}\},\\
\Vert \chi_{{\mathbb A},\infty} \Vert &=& \min \{ |{\cal F}| \such
{\cal F} \subseteq \ell^\infty \: \wedge \\
&&  \makebox[2cm]{}
(\forall A \in
{\mathbb A}) \;  (\exists f \in
{\cal F}) \; A \lim f \mbox{ does not exist}\}.
\end{eqnarray*}
By replacing $\ell^\infty$ by ${}^\omega 2$, the set of 
$\omega$-sequences with values in 2, we get the variations
$\chi_{{\mathbb A},2}$.
In \cite{Mi6} we showed that
for the cardinals we are interested in, ${}^\omega 2$ and $\ell^\infty$ 
give the same result. From now on we shall work with
${}^\omega 2$.

\smallskip

\Vo\ (cf.\ \cite{Vo2}) also gave some bounds valid for any ${\mathbb A}$
that contains at least all matrices which have exactly 
one non-zero entry in each line:
\[ 
\gs \leq \Vert\chi_{{\mathbb A},2}\Vert \leq \gb
\cdot \gs.\]
We write $\chi$ for $\Vert\chi_{{\mathbb M},2}\Vert$.

\smallskip

In \cite{Mi6} we showed that  $\chi < \gb \cdot \gs$
is consistent  relative to \zfc.
Here, we show the complementary consistency result, that
 $\gs < \chi$ is consistent. We get the convergence with positive matrices.

\smallskip
 
Now we recall here the definitions of the cardinal characteristics
$\gb$ and $\gs$
involved:
The order of eventual dominance $\leq^\ast$ is defined as follows: For
$f,g \in \omega^\omega$ we say $f \leq^\ast g$ if there is $k \in
\omega$ such that for all $n \geq k$ we have $f(n) \leq g(n)$.

The unbounding number $\gb$ is the smallest size of a subset ${\cal B}
\subseteq {}^\omega \omega$ 
such that for each $f \in {}^\omega \omega$ there is
some $b \in {\cal B}$ such that $b \not\leq^\ast f$. The splitting number
$\gs$ is the smallest size of a subset ${\cal S} \subseteq [\omega]^\omega$
such  that for each $X \in \potinf$ there is some $S \in {\cal S}$ such
that $X \cap S$ and $X \setminus S$ are both infinite. The latter is
expressed as ``$S$ splits $X$'', and $\cal S$ is called a splitting family.
For more information on these cardinal characteristics,
we refer the reader to the survey articles \cite{Blasshandbook,
vanDouwen, Vaughan}.

\smallskip

If $A \lim f$ exists, then also $A' \lim f$ exists for any $A'$ that
is gotten from $A$ be erasing rows and moving the remaining
(infinitely many) rows together. We may further change $A'$ by
keeping only finitely many non-zero entries in each row, such that the 
neglected ones have a negligible absolute sum, and 
then possibly multiplying the remaining ones such that they again sum up to 1.
Hence, after possibly further deleting
of lines  we may restrict the set of Toeplitz matrices 
to linear Toeplitz matrices. A matrix
is linear iff each column $j$ has at most one entry $a_{i,j} \neq 0$ and
for $j < j'$ the $i$ with $a_{i,j} \neq 0$ is smaller or equal to the
$i$ with $a_{i,j'} \neq 0$ if both exist, in picture
\begin{small}
\begin{equation*}
\begin{pmatrix}
c_0(0) & \dots c_0(\mup(c_0) -1) & 0 & \dots & 0 & 0 & \dots\\
0 & \dots            0           & c_1(\mdn(c_1)) & \dots & c_1(\mup(c_1) -1)
& 0 & \dots\\
0 & \dots 0 & 0 &\dots& 0 &  c_2(\mdn(c_2)) &\dots \\
\vdots 
\end{pmatrix}
\end{equation*}
\end{small}
\smallskip

Linear matrices can be naturally (as in the picture) 
 read as $(c_n \such n \in \omega)$
where 
$c_n \colon [\mdn(c_n), \mup(c_n))
\to [0,1]$, $c_n(j) = a_{n,j}$, give the finitely many 
non-zero entries in row $n$, and
$\mup(c_{n-1}) = \mdn(c_n)$. The $c_n$ are special instances
of the 
weak creatures in the sense of \cite{RoSh:470}.
In the next two sections we shall show:
The $c_n$'s coming from the trunks of the conditions in 
the generic filter of our forcing $Q$ 
 give matrices that  make, after multiplication, members of 
${}^\omega 2$ from the ground model and members of ${}^\omega 2$ of 
any random 
extension convergent.

\section{A creature forcing}\label{S1}

In this section, we give a self-contained description of the 
creature forcing $Q$ which is the main tool for building the two kinds
of models in the next section.
Moreover, we explain the connections
and give the references to \cite{RoSh:470},
so that the reader can identify it as a special case of
an extensive framework.

\begin{definition}\label{1.1} 
a) We define a notion of forcing $Q$. Its members $p$ are of the form
$p=
(n,c_0,c_1, \dots )= (n^p, c_0^p, c_1^p, \dots )$ such that
\begin{myrules}
\item[(1)] $n^p \in \omega$.
\item[(2)] For each $i \in \omega$ there are 
$\mdn(c_i) < \mup(c_i) < \omega$ such that\\
$c_i \colon [\mdn(c_i),\mup(c_i)) \to [0,1]$, such that
$(\forall k \in \dom(c_i))(c_i(k) \cdot k! \in {\mathbb Z})$.
\item[(3)] $w(c_i) = \{ k \in [\mdn(c_i),\mup(c_i)) \such c_i(k)
\neq 0 \}$, and $\sum_{k \in w(c_i)} c_i(k) = 1$.
\item We let $\norm(c_i) = \mdn(c_i)$.
We denote by $K$ the set of those $c_i$.
\item[(4)] $\mup(c_i) = \mdn(c_{i+1})$.
\end{myrules}

We let $p \leq q$ (``$q$ is stronger than $p$'', we follow the
Jerusalem convention) if
\begin{myrules}
\item[(5)] $n^p \leq n^q$.
\item[(6)] $c_0^p = c_0^q, \dots , c^p_{n^p-1} = c^q_{n^p-1}$.
\item[(7)] there are $n^p \leq k_{n^p} < k_{n^p+1} < \dots $ and there are 
non-empty sets $u \subseteq [k_n,k_{n+1})$ 
and rationals
$d_\ell >0$ for $\ell \in u$ such that
$c^q_n = \sum\{d_\ell \cdot c_\ell^p \such \ell \in u\}$
and $\sum_{\ell \in u } d_\ell = 1$.
We let $ \Sigma(\langle c_\ell \such \ell \in [k_n, k_{n+1})
\rangle$ denote the collection of all $c_n$ gotten with any
$u \subseteq [k_n,k_{n+1})$ and any weights 
$d_\ell$ for $\ell \in u$. Thus $\mdn(c^q_n) = \mdn(c^p_{k_n})$
and $\mup(c^q_n)= \mup(c^p_{k_{n+1} -1})$. 
\end{myrules}

b) We write $p \leq_i q$ iff
$n^p = n^q$ and $c_j^q = c_j^p$ for $j < n^p + i$.

\end{definition}

\begin{remark}
The notation we used in \ref{1.1} is natural to describe
our forcing in a compact manner. However, it does not coincide with
the notation given for the general framework in \cite{RoSh:470}. Here
is a translation:
We write \\ $((c_0^p,c_1^p, \dots c_{n-1}^p), c_n^p, c_{n+1}^p, \dots)$
instead of $(n^p, c_0^p, c_1^p, \dots)$, which contains the 
same information. Then we write
\begin{equation}\tag{$\ast$}\label{ast}
((c_0^p,c_1^p, \dots c_{n-1}^p), c_n^p, c_{n+1}^p, \dots)
= (w^p, t_0^p, t_1^p \dots).
\end{equation}
Then the $t_i^p$ are (simple cases) of components of weak creatures in the
sense of \cite[1.1.1 to 1.1.10]{RoSh:470}.
If we write $\bf t = (\norm({\bf t}), \val({\bf t}), \dis({\bf t}))$
for a weak creature in the sense of \cite{RoSh:470}, then we have 
that $\dis$ is the empty function, and 
$t_i$ is part of such a ${\bf t}$ in the following sense:
$\norm({\bf t}) = \mdn(t_i)$, $\rge(\val({\bf t})) = \{ t_i\}$.
We set ${\bf H}(i) = \left\{0,\frac{1}{i!},\frac{2}{i!},
\dots , \frac{i! - 1}{i!}, 1\right\}$ and $t_i \in \prod_{m \in 
[\mdn(t_i), \mup(t_i))} {\bf H}(i)$. $K$ is a collection 
of weak creatures, and $\Sigma$ from \ref{1.1}(7) 
is a composition operation. Thus
our $Q$ is $Q^*_{s\infty}(K,\Sigma)$ in Ros\l anowski's and Shelah's 
framework  and is finitary
and nice and satisfies some norm-conditions. 
We do not give the definitions of
these properties, because we are working with our specific case. 
The interested reader should consult \cite{RoSh:470}.
We use $w$, $w^p$, $w^q$ for the trunks in the representation as in 
\eqref{ast}.

\end{remark}

In order to make our work self-contained, we write a proof that
$Q$ allows continuous reading
of names and hence is proper. In this section,
we use the notation as in \eqref{ast},
because it is more suitable.

\begin{definition} \label{1.3}
$q=(w^q, t_0^q, \dots)$ approximates $\name{\tau}$ at $t_n^q$ iff
for all $r$ (if $q \leq r$ and $r$ forces a value to $\name{\tau}$, 
then $r^{q,n}$ forces this, where
$t_i^{r^{q,n}} = t_i^r$ for $i < n$ and
$\{ t_i^{r^{q,n}} \such i \geq n \} =
\{ t_i^q \such i < \omega, \mdn(t_i^q) \geq \mup(t^{r}_{n-1}) \}.$)
\end{definition}

\begin{definition}\label{1.4}
For $w\in\bigcup\limits_{m<\omega}\prod\limits_{i<m}{\bf H}(i)$ and
${\mathcal S}\in [K]^{\textstyle {\leq}\omega}$ we 
define the set $\pos(w,{\mathcal S})$ of
possible extensions of $w$ from the point of view of ${\mathcal S}$ (with respect to
$(K,\Sigma)$) as: 
\begin{eqnarray*}
\pos^*(w,{\mathcal S}) & =&
\Sigma({\mathcal S})\\
(&=& \{u: (\exists s\in\Sigma({\mathcal S}))(\langle
w,u\rangle\in\val[s])\}
\\
&&\mbox{for a general creature forcing}), 
\end{eqnarray*}
     
\[\hspace{-1cm}
\begin{array}{ll}
\pos(w,{\mathcal S})=\{u:&\!\!\!\!\mbox{there are disjoint sets }{\mathcal S}_i\mbox{ (for 
$i<m<\omega$) with }\bigcup\limits_{i<m}{\mathcal S}_i={\mathcal S}\\
\ &\mbox{and a sequence }0<\ell_0<\ldots<\ell_{m-1}<\lh(u)\mbox{ such that}\\
\ & u{\restriction} \ell_0\in\pos^*(w,{\mathcal S}_0)\ \&\ \\
\ & u{\restriction} \ell_1\in\pos^*(u{\restriction}\ell_0,{\mathcal S}_1)\ \&\ \ldots\
\&\ u\in\pos^*(u{\restriction} \ell_{m-1},{\mathcal S}_{m-1})\}.\\  
\end{array}
\]
\end{definition}

\begin{lemma}\label{1.5}
\label{deciding}
 (The case $\ell=0$ of \cite[Theorem 2.1.4]{RoSh:470}) 
$Q$ has continuous reading of names, i.e.
if $p \Vdash \name{\tau} \colon \omega \to V \mbox{ (old universe) }$
there is $q=(w^q,s_0, s_1 \dots )$  such that 
\begin{myrules}
\item[$(\alpha)$] $p \leq_0 q \in Q$,
\item[$(\beta)$] if $n < \omega$ 
and $m \leq \mup(s_{n-1})$ then the condition $q$ approximates
$\name{\tau}(m)$ at $s_n$.
\nothing{
d
 $q \leq r$ and $r$ forces a value to $\name{\tau}(i)$
and $i < \mup(t^r_{n_r -1})$, then $r^q$ forces this, where
$n^{r^q} = n^r$, $t_n^{r^q} = t_n^r$ for $n < n^r$ and
$\{ t_n^{r^q} \such n \geq n^r \} =
\{ t_n^q \such n < \omega, \mdn(t_n^q) > \mup(t^{r^q}_{n^{r^q}-1} \}$.
}
\end{myrules}
\end{lemma}

\proof
Let $p = (w^p,t_0^p,t_1^p, \dots )$. Let $w^q = w^p$.  
Now, by induction on $n \geq 0$ 
we define $q_n, s_n, t^n_{n+1}, t^n_{n+2}, \dots$ such that:
\begin{myrules}
\item[(i)] $q_0 = p$,
\item[(ii)] $q_{n+1} = ( w^p,s_0,\dots,s_n,t^n_{n+1},t^n_{n+2},\dots)
\in Q$,
\item[(iii)] $q_n \leq_n q_{n+1}$,
\item[(iv)] if $w_1 \in {\rm pos}(w^p,s_0,\dots s_{n-1})$, and $m \leq
\mup(s_{n-1})$ and there is a condition $r \in Q$, $r \geq_0
(w_1,s_n,t^n_{n+1},t^n_{n+2}, \dots )$ which decides the value of
$\name{\tau}(m)$ then the condition 
$(w_1,s_n,t^n_{n+1},t^n_{n+2}, \dots )$
already does it.
\end{myrules}

Arriving at stage $n \geq 0$ we have defined 
$$q_n=(w^p,s_0,s_1,\dots,s_{n-1},t^{n-1}_n,t^{n-1}_{n+1},\dots).$$
Let $\langle (w^n_i, m^n_i):
i< K_n\rangle$ be an enumeration of
\[{\rm pos}(w^p,s_0,\ldots,s_{n-1})\times(\mup(s_{n-1})+1)\]
(since each ${\bf H}(m)$ is finite, $K_n$ is finite). Next choose 
by induction on
$k\leq K_n$ conditions $q_{n,k}\in Q$ such that:
\begin{myrules}
\item[$(\alpha)$] $q_{n, 0}=q_n$. 
\item[$(\beta)$]  $q_{n,k}$ 
is of the form $(w^p,s_0,\dots,s_{n-1},t_n^{n,k},
t^{n, k}_{n+1}, t^{n, k}_{n + 2},\ldots)$.
We set $w^n_k = (w^p,s_0,
\dots s_{n-1})$. 
\item[$(\gamma)$] $q_{n,k}\leq_n q_{n,k+1}$.
\item[$(\delta)$] If, in $Q$, 
there is a condition
$r\geq_0 (w^n_k,t^{n,k}_n,t^{n,k}_{n+1},t^{n,k}_{n+2},\ldots)$ which decides
(in $Q$) the value of $\name{\tau}(m^n_k)$, then
\[(w^n_k,t^{n,k+1}_n,t^{n,k+1}_{n+1},t^{n,k+1}_{n+2},\ldots)\in Q\]
is a condition which forces a value to $\name{\tau}(m^n_k)$.    

\end{myrules}
For this part of the construction we need our standard assumption that
we may iterate the process in \ref{1.1}(7). Note, that choosing
$(w^n_k,t^{n,k+1}_n,t^{n,k+1}_{n+1},t^{n,k+1}_{n+2},\ldots)$ we want to be 
sure that 
\[(w^p,s_0,\ldots,s_{n-1},t^{n,k+1}_n,t^{n,k+1}_{n+1},t^{n,k+1}_{n+2},\ldots)
\in Q.\]
Next, the condition $q_{n+1}\stackrel{\rm def}{=}
q_{n,K_n}\in Q$ satisfies (iv): the keys are the
clause ($\delta$) and the fact that 
\[(w^n_k,t^{n,k+1}_n,t^{n,k+1}_{n+1},t^{n,k+1}_{n+2},\ldots)\leq
(w^n_k,t^{n,K_n}_n, t^{n,K_n}_{n+1},t^{n,K_n}_{n+2},\ldots)\in
Q.\]
Thus $s_n\stackrel{\rm def}{=} t^{n,K_n}_n$, 
$t^{n+1}_{n+k}\stackrel{\rm def}{=} t^{n,K_n}_{n+k}$ 
and
$q_{n+1}= (w^p,s_0, \dots, s_n,t^{n+1}_n, \dots )$
 are as required. 

Now, by a fusion argument
\[q\stackrel{\rm def}{=}(w^p,s_0,s_1,\ldots,s_l,s_{l+1},\ldots)=\lim_n q_n\in
Q.\] 
It is
easily seen that $q$ satisfies the assertions of the theorem. \proofend
\nothing{
\begin{definition}\label{1.6}
For $p,q \in Q$ we write $q \leq_{apr} q$ if there is 
some $n$, $w^q$ such that
$q=(w^q,t_n^p,t_{n+1}^p,\dots)$ and $p \leq q$. 
That is, $q$ is only in the
trunk stronger than $p$.
\end{definition}
}

\begin{lemma}(\cite[Corollary 2.1.6]{RoSh:470})
\begin{myrules}
\item[(a)] Suppose that $\name{\tau}_n$ are $Q$-names 
for ordinals and $q\in Q$ is a
condition satisfying $(\beta)$ of \ref{deciding}. Further
assume that $q\leq r\in Q$ and
$r\Vdash $``$\name{\tau}_m=\alpha$'' (for some ordinal $\alpha$).\\ 
Then $q'= r^{q,m}$ forces this.

\item[(b)] The forcing notion $Q$ is proper.
\end{myrules}
\end{lemma}

\proof
(a) is a special case of the previous lemma.

For (b), we use the  equivalent definition of properness given in 
\cite[III.2.13]{Sh:h}, and the fact that
$\{q'\!\in Q: (\exists r \geq q)(\exists n)
q'=r^{q,n}\}$ 
is countable provided $\bigcup\limits_{i<\omega}{\bf H}(i)$ 
is countable.

\section{The effect of $Q$ on random reals}\label{S2}

Let $G$ be $Q$-generic over $V$.
We set $c^G_n = c_n^q$ for $q \in G$ and $n^q >n$. This is well defined.
  Let $\name{c_n}$  be a name for it.
Our aim is to show that 
multiplication by the matrix whose $n$-th row is $c_n$ 
makes any real from the ground model and even any real
from a random extension of the ground model  convergent.
For background information about random reals we refer the reader to 
\cite[\S 42]{Jech}. The Lebesgue measure ist denoted by $\Leb$.
With ``adding $\kappa$ random reals'' we mean forcing with
the measure algebra $R_\kappa$ on $2^{\omega \times\kappa}$, that is adding 
$\kappa$ random reals at once or ``side-by-side'' and not successively.

\begin{definition}\label{2.1}
\begin{myrules}
\item[(1)] Let $\may_k(p) = \{ c_n^r \such p \leq_k r, n \geq n^p + k\}$.
\item[(2)] For a creature $c$ and $\eta \in {}^\omega 2$ let $\aver(\eta,c) 
= \sum_{k \in w(c)} c(k) \eta(k)$.
\end{myrules}
\end{definition}

\begin{mainlemma}\label{2.2}
Assume that 
\begin{myrules}
\item[(A)] $\name{\eta}$ is a random name of a member of ${}^\omega 2$,
$ \name{\eta} = 
f(\name{r})$ where $f$ is Borel and $\name{r}$ is as name of the 
random generic real,
\item[(B)] $p \in Q$,
\item[(C)] $k^* < \omega$.
\end{myrules}
Then for every $k \geq k^*$ there is some $q(k) \in Q$ such that
\begin{myrules}
\item[($\alpha$)]$p \leq_{k^*} q(k)$,
\item[$(\beta)$] for all $\ell$,
if $k^* \leq k < \ell <  \omega$ and 
$c_1, c_2 \in \may_\ell(q(k))$ then 
  $$ \frac{1}{\ell!} >
\Leb\left\{r \such \frac{3}{2^k} \leq \left| \aver(f(r), c_1) - 
\aver(f(r), c_0) \right| \right\}.$$

\end{myrules}
\end{mainlemma}

\proof
For $q \in Q$ and $k, \ell \in \omega$, 
$i \in \{0,1, \dots, 2^k\}$ we set
\begin{eqnarray*}
\err_{k,i}(\name{\eta},c) &=& \Expect\left(\left|\aver(\name{\eta},c) - 
\frac{i}{2^k}\right|\right)\\ & = & \int_0^1
\left|\aver(f(r),c)
- \frac{i}{2^k} \right| \; d \Leb(r),\\
\eee^\ell_{k,i}(\name{\eta},q) &=& \inf \{ \err_{k,i}(\name{\eta},c)
\such c \in \may_\ell(q)\}.
\end{eqnarray*} 

Note that $\err_{k,i}(\name{\eta},c)$ is a real and no longer a random name.
So the infimum is well-defined.

Now, if $\ell_1 < \ell_2$ then $\may_{\ell_1}(q) 
\supseteq \may_{\ell_2}(q)$ and hence
\begin{equation*}\label{mono}
\eee^{\ell_1}_{k,i}(\name{\eta},q) \leq \eee^{\ell_2}_{k,i}(\name{\eta},q).
\end{equation*}

So $\langle \eee^\ell_{k,i}(\name{\eta},q) \such \ell \in \omega
\rangle$ is an
increasing bounded sequence and 

\begin{equation*}
\eee^*_{k,i}(\name{\eta},q) =
\lim \langle \eee^\ell_{k,i}(\name{\eta},q) \such \ell \in \omega
\rangle
\end{equation*}
 is well-defined.

We fix $i \leq 2^k$, until Subclaim 4, when we start 
looking at all $i$ together.
\nothing{, and later we shall vary $k$ as well.}

\smallskip

Subclaim 1: 
There is some $q^{k,i}_1=q_1 \geq_{k^\ast} p$ such that
for $\ell \geq k^\ast$
$$\eee^*_{k,i}(\name{\eta}, p) - \frac{1}{\ell} 
\leq \err_{k,i}(\name{\eta}, c^{q_1}_\ell)
\leq \eee^*_{k,i}(\name{\eta}, p) + \frac{1}{\ell}.$$
Moreover, if $\mdn(c_{\ell'}^{q_1}) = \mdn(c^p_{\ell})$ then
$\eee^{\ell'}_{k,i}(\name{\eta},q_1) \geq
\eee_{k,i}^*(\name{\eta},p)-\frac{1}{\ell}$.

Why? We choose $c_\ell^{q_1}$ by induction 
on $\ell$: For $\ell \leq n^p + k^\ast$, we take
$c_\ell^{q_1} = c_\ell^p$.
Suppose that we have chosen $c_m^{q_1}$ for $m < \ell$ and that
we are to choose $c_\ell^{q_1}$, $\ell > n^p + k^\ast$. 
We set $\eps = \frac{1}{\ell}$.
By possibly end-extending $c_{\ell -1}^{q_1}$ by zeroes
we may assume
that $\mup(c_{\ell -1}^{q_1}) = \mup(c_{\ell'}^p)$ for such a large $\ell'
\geq \ell$ such that for all $\ell'' \geq \ell'$,
$\eee_{k,i}^{\ell''}(\name{\eta},p) \geq \eee_{k,i}^\ast(\name{\eta},p)
- \eps$. Then we take $c_\ell = c_\ell^{q_1} \in \may_{\ell''}(p)$ 
such that $\err_{k,i}(\name{\eta},c_\ell^{q_1}) 
\leq \eee_{k,i}^{\ell''}(\name{\eta},p)  + \eps
\leq \eee_{k,i}^\ast(\name{\eta},p) + \eps$. 
On the other side we have that 
$\err_{k,i}(\name{\eta},c_\ell^{q_1}) 
\geq \eee_{k,i}^{\ell''}(\name{\eta},p) 
\geq \eee_{k,i}^\ast(\name{\eta},p) - \eps$. The fact that
 this holds also for  $\ell' \leq \ell$ if  $\mdn(c_{\ell'}^{q_1}) = 
\mdn(c^p_{\ell})$ yields the ``moreover'' part.
\medskip

Subclaim 2: In Claim 1, if $\ell \geq k^*$ and $q^{k,i}_1 \leq_\ell q_2$ then
$$\eee^*_{k,i}(\name{\eta},q) - \frac{1}{\ell}
\leq \err_{k,i}(\name{\eta}, c_\ell^{q_2}) \leq
\eee^*_{k,i}(\name{\eta},q) + \frac{1}{\ell}.$$.

Why? By the definition if suffices to show: 
\begin{equation}\tag{$\otimes$}\label{otimes}
\begin{split}
& \mbox{if } \ell_1 < \cdots < \ell_t < \omega
\mbox{ and } d_1 , \dots d_t \geq 0 \mbox{ and } d_1 + \cdots + d_t =1,\\
& \mbox{and } c_\ell^{q_2}
=d_1 c_1^{q_1} + \cdots + d_t c_t^{q_1},\\
&\mbox{then } 
\eee^*_{k,i}(\name{\eta},q_1) - \frac{1}{\ell}
\leq \err_{k,i}(\name{\eta}, c_\ell^{q_2}) \leq
\eee^*_{k,i}(\name{\eta},q_1) + \frac{1}{\ell}.
\end{split}
\end{equation}

The first inequality holds by the ``moreover'' after the first inequality in 
the previous claim. For the second inequality it suffices to show that
$$\err_{k,i}(\name{\eta}, c) \leq \sum_{s=1}^t
d_s \err_{k,i}(\name{\eta}, c_s^{q_1}).$$
For this is suffices to show that 
$$\Expect\left(\left|\aver(\name{\eta},c) - \frac{i}{2^k} \right|\right) 
\leq 
\sum_{s=1}^t d_s
\Expect\left(\left|\aver(\name{\eta},c_{s}^{q_1}) - 
\frac{i}{2^k} \right|\right), 
$$
and writing this explicitly noting 
that $\Expect$ is actually a Lebesgue integral and that $d_s \geq 0$ and 
that $\sum_{s} d_s = 1$ we finish by the triangular 
inequality.

\medskip

Subclaim 3: Let $q^{k,i}$ be as in Subclaim 2.
For all $\ell$, if $c_0, c_1 \in \may_\ell(q^{k,i}_1)$, then
$$ \frac{2^{k+1}}{\ell} \geq
\Leb\left\{ r \such \aver(f(r), c_0) \geq \frac{i+1}{2^k} 
\wedge \aver(f(r),c_1) \leq \frac{i-1}{2^k} \right\}.$$

Why? Consider $c = \frac{1}{2} c_0 + 
\frac{1}{2} c_1 \in \may_\ell(q_1)$ .
Write \\
$A=\left\{ r \such \aver(f(r), c_0)
 \geq \frac{i+1}{2^k} \wedge \aver(f(r),c_1) \leq \frac{i-1}{2^k} \right\}$.

\begin{eqnarray*}
\frac{2}{\ell} & \geq & \frac{1}{2} \err_{k,i}(\name{\eta},c_0)
+  \frac{1}{2} \err_{k,i}(\name{\eta},c_1) -\err_{k,i}(\name{\eta},c)\\
& = & 
\int_0^1
\left(\frac{1}{2}\left|\aver(f(r),c_0)
- \frac{i}{2^k} \right|   +\right.
\frac{1}{2}\left|\aver(f(r),c_1)
- \frac{i}{2^k} \right| -\\
&&
\left.\left|\aver(f(r)),c)
- \frac{i}{2^k} \right|\right) d \Leb(r)\\
&\geq& 
\int_A
\left(\frac{1}{2}\left|\aver(f(r),c_0)
- \frac{i}{2^k} \right|   +\right.
\frac{1}{2}\left|\aver(f(r),c_1)
- \frac{i}{2^k} \right| -\\
&&
\left.\left|\aver(f(r)),c)
- \frac{i}{2^k} \right|\right) d \Leb(r)\\
&\geq& 
\frac{1}{2^{k}} \Leb(A).
\end{eqnarray*}

\nothing{
We divide the truth value in the defintion of $\err_{k,i}(\name{\eta},c)$ by the events above.
Outside it is $\leq \frac{1}{2} \err_{k,i}(\name{\eta}, c_1) + \frac{1}{2} \err_{k,i}(\name{\eta},c_2)
\leq \err^*_{k,i} + \frac{1}{\ell}$.
Inside the result drops by $\Leb$(event above) $\times \frac{1}{2^{k-1}}$
as $d_1 \leq -\frac{1}{2^k}$ and $d_2 \geq \frac{1}{2^k}$ implies that $d_1 + d_2| \leq |d_1| +|d_2| - \frac{1}{2^{k_1}}$.
}

\medskip

Subclaim 4: For every $q \in Q$ and $k^*$ we can find $q^{k}$ such that
\begin{myrules}
\item[$\alpha$)] $q \leq_{k^*} q^{k}$,

\item[$\beta$)]  if $\ell \in [k,\omega)$ and $c_0,c_1 \in 
\may_\ell(q^{k})$ and $i \in \{1,2, \dots, 2^k-1\}$ then 
$\frac{2^{k+1}}{\ell} 
> \Leb\left\{ r \such \aver(f(r), c_0) 
\geq \frac{i+1}{2^k} \wedge \aver(f(r),c_1) \leq \frac{i-1}{2^k} \right\}.$

\item[$\gamma$)] This holds also for every $q^{*} \geq q^{k}$.
\end{myrules}
Why? Repeat Subclaims 1 and 2 and 3 choosing $q^{k,i}$, $i = 0,1, 
\dots, 2^k$. We let  $q_0 = q$
and choose $q^{k,i+1} $ such that it
 relates to $q^{k,i}$ like $q_1$ to $q$.
 
Now $q^{k}=q^{k,2^k}$ is o.k. 
Note that according to \eqref{otimes} 
thinning and averaging can only help.

\medskip

Subclaim 5: Let $q^k$ be as in Subclaim 4.
For $\ell \geq k$ there
is  $q(k,\ell)\geq_{\ell-1}q^k$ such that  
for $c_0,c_1 \in \may_\ell(q(k,\ell))$,
$$ \frac{1}{\ell!} >
\Leb\left\{r \such \frac{3}{2^k} \leq \left| \aver(f(r), c_1) - 
\aver(f(r), c_0) \right| \right\}.$$

Why? The event $\frac{3}{2^k} \leq \left| \aver(f(r), c_1) - 
\aver(f(r), c_0) \right|$  implies 
that for some $i \in \{1,2,\dots,2^k-1\}$ we have 
$\aver(f(r), c_1) \geq \frac{i+1}{2^k}
\wedge \aver(f(r), c_2) \leq \frac{i-1}{2^k}$ or
vice versa. So it is incuded in the union of 
$2 \times (2^k -1)$ events, each of measure $\leq \frac{2^{k+1}}{\ell}$.
 Hence it itself has measure
$\leq \frac{2^{2k+2}}{\ell}$. By thinning out 
$q^k$ (by  moving the former $\ell$ far out
by putting in a lot of zeroes and thus 
having as new $c_\ell$'s weak creatures that
were formerly labelled with a much larger $\ell$ and thus giving 
a much smaller quotient according to Subclaim 4) 
we replace 
$\frac{2^{2k+2}}{\ell}$ by $\frac{1}{\ell!}$.

\medskip

Subclaim 6: Finally we come to the
$q(k)$ from part $(\beta)$ of the lemma:
 For any $k$  there is $q(k)$ such that
$q \leq_{k^*} q(k)$ and for any $\ell\geq k$ and any 
 $c_1, c_2 \in \may_\ell(q^*)$ then 
 
  $$ \frac{1}{\ell!} >
\Leb\left\{r \such \frac{3}{2^k} \leq \left| \aver(f(r), c_1) - 
\aver(f(r), c_0) \right| \right\}.$$

Why? Like in the previous claim we choose inductively $q(k,\ell)$ such that 
$q_0 =p$ and $q(k,\ell+1) \geq_\ell q(k,\ell)$ 
and $(q(k,\ell+1),q(k,\ell), \ell)$ are like 
$(q(k,\ell), q, \ell)$ from Subclaim 5, but for larger and larger $\ell$. Now
$$q(k) =(n^p + k, c_0^p, \dots, c_{n^p + k}^p,
c_{n^p + k +1}^{q(k,n^p+k+1)}, c_{n^p +k+1}^{q(k,n^p+k+2)}, \dots) $$
is as required in $(\alpha)$ and $(\beta)$
of the conclusion; we have even $q(k) \geq_k p$.
\proofend

\begin{conclusion}\label{2.3}
$\Vdash_Q $ ``if $\name{\eta} \in V$ is a random name of a
 member in $2^\omega$ (i.e.\ a name for a real in $V^{R_\omega}$) 
then ``$\Vdash_{R_\omega} \langle \aver(\name{\eta}, 
\name{c}_n) \such n \in \omega \rangle$ converges'' ''
\end{conclusion}

\proof 
Let $q \in Q$  and $\eps > 0$ be given. Let $\name{\eta} = f(\name{r})$,
$f \in V$,
be a random name for a real. We take 
$k_0$ such that $\frac{3}{2^k} < \eps$.
Then we take for $q(k)\geq q$ as in the Main Lemma.
We set
$$A_{k,c_0,c_1} = \left\{ r  \such 
\frac{3}{2^k} > \left| \aver(f(r), c_1) - 
\aver(f(r), c_0) \right| \right\}.$$
Since $\sum_{\ell \geq 1} \frac{1}{\ell !} < \infty$, we 
can apply the Borell Cantelli lemma and get:

For any sequence $\langle c_\ell \such \ell \in \omega \rangle$
such that $c_\ell \in \may_\ell(q^*(k))$  we have that
$$\Leb \left(\bigcup_{K \in [k,\omega)} \bigcap_{\ell \geq K} 
A_{k_0,c_\ell,c_{\ell+1}} \right)
=1.$$
So $r \in \bigcap_{\ell \geq K} A_{k,c_\ell,c_{\ell+1}}$ for some 
$K \geq k$.
So $q(k)$ forces that $\langle c_\ell \such \ell \in \omega \rangle$ 
describes a matrix whose product with $\eta$ lies eventually
within an $\eps$ interval.
Now we take smaller and smaller $\eps$'s and a density argument.
\proofend

\begin{conclusion} \label{2.4}
Let $P_{\omega_2}=\langle P_i, \name{Q}_j \such i \leq \omega_2, 
j < \omega_2 \rangle$ be a countable support iteration of 
$\name{Q}_i$, where $Q_i$ is $Q$ defined in $V^{P_i}$, and 
let $\name{R}_{\omega_1}$ be a 
$P_{\omega_2}$ name of the $\aleph_1$-random algebra.
Then in $V^{P_{\omega_2} \ast \name{R}_{\omega_1}}$ we have 
$ \gs = \aleph_1$ and $\chi > \aleph_1$.
\end{conclusion}

\proof Dow proves in \cite[Lemma 2.3]{Dow} that
 $s = \aleph_1$ after adding $\aleph_1$ or more random reals,
over any ground model. In order to
show $\chi > \aleph_1$, let $\eta_i$ , $i < \omega_1$ be reals in
$V^{P_{\omega_2} \ast \name{R}_{\omega_1}}$. Over $V^{P_{\omega_2}}$,
each $\eta_i$ has a $R_{\omega_1}$-name $\name{\eta_i}$.
 Since the random algebra is c.c.c, there are w.l.o.g.\ only countably
many of the $\aleph_1$ random reals mentioned in $\name{\eta_i}$.
Let $\name{\eta'_i}$ be got from $\name{\eta_i}$ by replacing these
countably many by the first $\omega$ ones and then doing as if
it were just one random real. This is possible because 
$R_1$ and $R_\omega$ are equivalent forcings.
Since the random algebra is c.c.c.,
the name $\name{\eta'_i}$ can be coded as a single real $r_i$
 in $V^{P_{\omega_2}}$.
Now, by \cite[V.4.4.]{Properforcing} and by the properness of the $
\name{Q_j}$, this name $r_i$ appears at some stage $\alpha(\eta_i)
<\aleph_2$ in the iteration $P_{\omega_2}$. We take the
supremum $\alpha$ of all the $\alpha(\eta_i)$, $i < \omega_1$. 
We apply the Main 
Lemma to the $\name{\eta'_i}$.  
Thus $Q_\alpha$ adds a Toeplitz matrix, that makes after multiplication all
the $\eta'_i$  convergent.
Since the Main Lemma applies to all random
algebras simultaneously, this matrix makes also the
$\eta_i$ convergent. 
\proofend

\begin{definition}
\label{2.5}
\begin{myrules}
\item[(1)]
 $Q_{pr}= \{p \in Q \such n^p = 0\}$ is
called  the pure part of $Q$.
\item[(2)]
We write $p \leq^* q$ if there are some $w$, $n$  such that
$p \leq (w,t_n^q,t_{n+1}^q \dots )$. So, it is up to a finite ``mistake''
$p \leq q$.
\end{myrules}
\end{definition}

\begin{fact}\label{2.6}
If $\langle p_i \such  i < \gamma \rangle$ is $\leq^*$-increasing 
in $Q$ and 
$\MA_{|\gamma|}$ holds, then there is $p \in Q_{pr}$ such that for 
all $i < \delta$,
$p_i \leq^*  p$.
\end{fact}
\proof
We apply $\MA_{|\gamma|}$ to the following partial order
$P$: Conditions are $(s,F)$ where
$s=( t_0^p, \dots, t_n^p)$ is an initial segment of a
condition in $Q_{pr}$ and $F \subset \gamma$ is a finite set.
We let $(s,F) \leq_P (t,G)$ iff
$s \trianglelefteq t$ and $F \subseteq G$ and  
$(\forall n \in \lgg(t) - \lgg(s))(\forall \alpha \in F)
(n > $ (all mistakes between the $p_\alpha)
\rightarrow t_n \in \Sigma(c_i^{p_\alpha} \such i \in {\mathcal S}(\alpha,n)$
for suitable ${\mathcal S}(\alpha,n)))$.
This forcing is c.c.c., because conditions with the same first component
are compatible and because there are only countably many possibilities
for the first component. It is easy to see that
for $\alpha < \delta$ the sets $D_\alpha = \{ (s,F) \such
\alpha \in F \}$ is dense and that for $n \in \omega$ the sets
$D^n =\{ (s,F) \such \lgg(s) \geq n \}$ are dense.
Hence if $G$ is generic, then
$p= \bigcup\{s \such \exists F (s,F) \in G \}
\geq^* p_\alpha$ for all $\alpha$.
\proofend

\begin{conclusion}\label{2.7}
If $V \models \MA_\kappa$ and $\kappa > \delta > \aleph_0$,
then in $V^{R_\delta}$ then matrix number is $\geq \kappa$ and the 
splitting number is $\aleph_1$.
\end{conclusion}

\proof As mentioned, \cite{Dow} shows the
the result on the splitting number.
For the matrix number, let random names
 $\name{\eta_i}$, $i
< \gamma$ be given in $V$, $\gamma < \kappa$.
 We fix $\eps > 0$ and $K$ as in
the proof of \ref{2.3}.
 We choose for $i < \gamma$, 
$p^i = \langle c^i_k \such k \in \omega \rangle$
 as in the end of the proof of
\ref{2.3} for $\name{\eta_i}$  and use 
and Fact 2.6. $\gamma+1$ times iteratively
and find a pure condition $p = \langle c_k \such k \in \omega\rangle
\geq^\ast p^i$ for all $i < \gamma$, 
that gives the lines of a matrix which brings everything into an
$\eps$-range. We denote these $c_k$ by 
 $c_k = c_k(\eps)$. Now by induction we choose $c_k$:
$c_0 = c_0(1)$, and 
$c_k= c_{k'}(\frac{1}{k'+1}) $ if $k' > k$ is the first $k''$ such
that $\mdn(c_{k''}(\frac{1}{k''+1})) > \mdn(c_{k-1})$.
The matrix with $c_k$ in the $k$th line
 acts as desired. (Now $\mup(c_k) > \mdn(c_{k+1})$ is possible
but this does not do any harm.)
\proofend

\smallskip

{\bf Acknowledgement:} The first author would like to thank Andreas Blass for
discussions on the subject and for reading and commenting.


\end{document}